\documentclass[12pt]{article}
\usepackage{amsgen, amsmath, amsfonts, amsthm, amssymb, enumerate,amscd,hyperref,graphicx, lipsum}
\usepackage[dvipsnames]{xcolor}
\hypersetup{colorlinks=true,linkcolor=blue,citecolor=magenta}
\usepackage[all]{xy}

\newtheorem{defn}{Definition}[section]
\newtheorem{prop}[defn]{Proposition}

\numberwithin{equation}{section}

\definecolor{blue}{rgb}{0,0,1}
\definecolor{red}{rgb}{1,0,0}
\definecolor{green}{rgb}{0,.6,.2}
\definecolor{purple}{rgb}{1,0,1}

\long\def\red#1\endred{\textcolor{red}{#1}}
\long\def\blue#1\endblue{\textcolor{blue}{#1}}
\long\def\purple#1\endpurple{\textcolor{purple}{ #1}}
\long\def\green#1\endgreen{\textcolor{green}{#1}}

\usepackage[OT2,T1]{fontenc}
\DeclareSymbolFont{cyrletters}{OT2}{wncyr}{m}{n}
\DeclareMathSymbol{\Sha}{\mathalpha}{cyrletters}{"58}

\newtheorem{lem}[defn]{Lemma}
\newtheorem{lemma}[defn]{Lemma}
\newtheorem{thm}[defn]{Theorem}
\newtheorem{theorem}[defn]{Theorem}
\newtheorem{conjecture}[defn]{Conjecture}
\newtheorem*{conj}{Conjecture}
\newtheorem*{conj*}{Conjecture}

\newtheorem{rem}[defn]{Remark}
\newtheorem{remark}[defn]{Remark}

\newcommand {\ZZ}{{\mathbb Z}}

\newcommand {\C}{{\mathbb C}}

\newcommand {\Q}{{\mathbb Q}}
\newcommand {\R}{{\mathbb R}}

\newcommand {\HH}{{\mathfrak  H}}

\def\dim{\operatorname{dim}}

\title{}
\author{ \\
}

\newcommand\blfootnote[1]{%
  \begingroup
  \renewcommand\thefootnote{}\footnote{#1}%
  \addtocounter{footnote}{-1}%
  \endgroup
}

\begin{document}

\title{Period polynomials, derivatives of $L$-functions, and zeros of polynomials }
\author{Nikolaos Diamantis (University of Nottingham)
\\
Larry Rolen (Hamilton Mathematics Institute \& Trinity College Dublin)
}

\maketitle

\centerline{{\it Dedicated to Don Zagier in honor of his 65th birthday}}

\begin{abstract}
Period polynomials have long been fruitful tools for the study of values of $L$-functions in the context of 
major outstanding conjectures. In this paper, we survey some facets of this study from the perspective of Eichler cohomology. We
discuss ways to incorporate non-cuspidal modular forms and values of derivatives of $L$-functions into the same framework.
We further review investigations of the location of zeros of the period polynomial as well as of its analogue for $L$-derivatives.
\end{abstract}

\section{Introduction}
The \blfootnote{\hskip -6mm 2010 \it Mathematics Subject Classification: \rm 11F11, 11F67.
\\
 \it Key words and phrases: \rm Periods of modular forms, Derivatives of $L$-functions, Eichler cohomology, Critical $L$-values.
}
period polynomial provides a way of encoding critical values of $L$-functions
associated to modular cusp forms that has proven very successful in the uncovering of
important arithmetic properties of $L$-values. As such, its structure and
properties as an object in its own right have attracted a lot of interest
from various
perspectives, one of the most important ones being that of
Zagier, as will become apparent below. To give an idea of the uses of the period polynomial and
its structure, we start by outlining its definition.

Let $f$ be an element of the space $S_k$ of weight $k$ cusp forms for SL$_2(\ZZ)$.
The {\it period polynomial} of $f$ is the polynomial in $X$ is given by
$$r_f(X):=\int_0^{\infty} f(\tau) (\tau-X)^{k-2} d\tau.$$
A relation with the $L$-function of $f$ is provided by the identity (cf eg. \cite{Z})
\begin{equation}
\label{LfPeriodP}
r_f(X)=\sum_{n=0}^{k-2} \binom{k-2}{n} i^{n-1}\Lambda_f(n+1)X^{n},
\end{equation}
where
$\Lambda_f(s):=(2 \pi)^{-s}\Gamma(s) L_f(s)$ is the ``completed'' $L$-function of $f.$

An example of the manner by which the structure of the period polynomial
leads to important arithmetic information about values of $L$-functions
is Manin's Periods Theorem.
The algebraic properties of $r_f$
({\it cocycle relations}) combined with the arithmetic nature of $f$
(as a Hecke eigenform) lead to a certain 1-dimensionality statement for $r_f$,
which, with \eqref{LfPeriodP} translates to the following proportionality
relation.
\begin{thm} Manin's Periods Theorem (\cite{M}).
Let $f$ be a normalized Hecke eigenform in $S_k$ with rational Fourier coefficients.
Then there exist $\omega_+(f)$, $\omega_-(f) \in \R$ such that
$$
\Lambda_f(s)/\omega_+(f), \quad \Lambda_f(w)/\omega_-(f) \in \mathbb{Q}
$$
for all $s,w$ with $1 \le s, w \le k-1$ and $s$ even, $w$ odd.
\end{thm} 
\begin{remark}
{\rm For illustration purposes, here we gave a special case of the actual theorem
which will be discussed in slightly more detail in the next section.}
\end{remark}

Fundamental applications such as the above theorem
have motivated closer independent
study of $r_f$, for instance as a polynomial.
The aspect we will be focusing
on in this survey is the location of the zeroes of $r_f(X).$

The strength of techniques based on $r_f$
has likewise motivated the search for
analogues of the period polynomial in other situations. The example
we will more closely be reviewing here is an analogue of the period polynomial
associated to {\it derivatives} of $L$-functions.

Derivatives of $L$-functions are the subject of some of the main
current conjectures in number theory, e.g. by Birch--Swinnerton-Dyer and by Beilinson.
To review a part of the latter, in an explicit formulation due to Kontsevich-Zagier
\cite{KZ}, we recall the definition of \emph{periods}, again in the form given in \cite{KZ}:
These are complex numbers whose real and imaginary parts have the form
$$\int_V \frac{P(\mathbf{x})}{Q(\mathbf{x})}d\mathbf{x},$$
where $V$ is a domain in $\mathbb{R}^n$ defined by polynomial inequalities with coefficients in $\Q$
and $P, Q \in \Q [X_1, \dots, X_n].$ The set $\mathcal{P}$ of periods contains $\pi,$
$\log(n)$ ($n \in \mathbb{N}$), etc. The arrangement of the following
special case of Beilinson's conjecture follows \cite{KZ}.
\noindent
\begin{conj}[Deligne-Beilinson-Scholl] \label{Beil}
Let $f$ be a weight $k$ Hecke eigencuspform for SL$_2(\mathbb{Z})$, $L_f(s)$ its $L$-function, and $m$ an integer.
Then, if $r$ is the order of vanishing of $L_f(s)$ at $s=m$, 
$$L^{(r)}(m) \in \mathcal P [1/\pi ].$$
\end{conj}
Apart from the cases $r=0$ (treated by Manin, Deligne, Beilinson, Deninger-Scholl; see \cite{KZ} and the references contained there)
and $r=1$ (thanks to Gross-Zagier \cite{GZ} and others), this conjecture is still open.
Analogues of the period polynomial for first derivatives of $L$-functions have been given in \cite{G, D}.
The version we will be using is
$$\int_0^{\infty} f(w)(w-z)^{k-2}\left(\log(w)-\frac{\pi i}{2}\right)dw.
$$
The justification for this choice will come from cohomological considerations (see Section \ref{Der}),
but, for the time being, we note that that this polynomial, in analogy with \eqref{LfPeriodP}, equals 
$$-\sum_{n=0}^{k-2} \binom{k-2}{n} i^{1-n}\Lambda_f'(n+1)z^{k-2-n}.$$

This polynomial has an algebraic structure
that fits into the same context as that of the standard period polynomial. It was recently observed by
the authors that, at least conjecturally, its zeros follow the same pattern as
those of the standard period polynomial.

\begin{conjecture}\label{Conj0} (``Riemann hypothesis for period polynomials attached to $L$-derivatives'') \cite{DR}
For any Hecke eigenform of weight $k$ on $\operatorname{SL}_2(\mathbb Z)$, and for each $m\in\mathbb Z_{\geq0}$, the polynomial
\[
Q_f(z):=\sum_{n=0}^{k-2} \binom{k-2}{n} i^{1-n}\Lambda_f^{(m)}(n+1)z^{k-2-n}
\]
has all its zeros on the unit circle.
Moreover, its odd part
\[\sum_{\substack{n=1 \\ n \, \, \text{odd}}}
^{k-3}\binom{k-2}{n} i^{1-n}\Lambda_f^{(m)}(n+1)z^{k-2-n}
\]
has all of its zeros on the unit circle, except for $0$, $\pm a,\pm1/a$ for some $a\in\mathbb R$. 
\end{conjecture}
 In \cite{DR}, this statement was proved in the case of Eisenstein series. 
 
In this survey, we will review the theory of period polynomials and of the ``period polynomials'' attached to
$L$-derivatives from a cohomological perspective. We will further survey conjectures and results about
zeros of period polynomials and of their counterparts for $L$-derivatives.

\section{Period polynomials}

\subsection{Period polynomials of cusp forms}

Set $\Gamma:=$PSL$_2(\mathbb{Z})$. This group is generated by $S=\left ( \begin{smallmatrix}
0 & -1 \\ 1 & 0 \end{smallmatrix} \right )$
and $T=\left ( \begin{smallmatrix} 1 & 1 \\ 0 & 1 \end{smallmatrix} \right )$ (or, more precisely, by
their images under the natural projection of SL$_2(\mathbb{Z})$ onto $\Gamma$). The only relations
are 
\begin{equation}
\label{ST}
S^2=(ST)^3=1.
\end{equation}
For $ \tau \in \mathfrak{H}$, let 
$f(\tau)=\sum_{n=1}^{\infty} a_n e^{2 \pi i n} $
be a cusp form of even weight $k$ for $\Gamma$.  A way to define the
period polynomial associated to $f$ is as a polynomial in $z$ of degree
$\le k-2$ given by
\begin{equation}\label{ppoly}
r_f(z):=\int_{0}^{\infty} f(\tau)(\tau-z)^{k-2}d\tau.
\end{equation}
The origin of this definition goes back (at least) to Poincar\'e (cf. \cite{P}) in the context of
work on abelian integrals (a fact brought to our attention by \cite{DIT}). Since then,
the period polynomial has been interpreted in several ways, each providing new insight and
leading to important applications. 
We will
review two that are most relevant for our purposes.

{\bf a. Eichler cohomology } \\
Firstly, Eichler \cite{EichlerVer} and Shimura \cite{ShimuraSur} viewed them
as periods of iterated integrals that are now called {\it Eichler integrals}:
\begin{equation}\label{EInt}
F(z)=(k-2)!\int_{\infty}^z \int_{\infty}^{z_1} \dots \int_{\infty}^{z_{k-2}} f(z_{k-1}) dz_{k-1} \dots d z_1=\int_{\infty}^z f(\tau) (\tau-z)^{k-2} d\tau.
\end{equation}
The relation of $F$ with $r_f$ is given by
\begin{equation}\label{cobound}
F(-1/z)z^{k-2}-F(z)=r_f(z).
\end{equation}
This identity can be viewed as the starting point of an algebraic approach
to the study of the period polynomial which has far-reaching implications. It
first implies that $r_f$ induces a $1$-cocycle in {\it Eichler cohomology}, which we will
now define. Since we will later need cocycles in more general cases, we recall the general
definition of cocycles.

We denote the space of $i$-cochains for $\Gamma$ with coefficients in a right $\Gamma$-module $M$
by $C^i(\Gamma, M).$ The differential $d^i\colon C^i(\Gamma, M) \to
C^{i+1}(\Gamma, M)$ is given by
\begin{align}\label{differential}
\begin{split}
&(d^i \sigma)(g_1, \dots, g_{i+1}):=\\
&
\sigma(g_2,\dots,g_{i+1}).g_1
+\sum_{j=1}^i (-1)^j \sigma(g_1,\dots, g_{j+1}g_j, \dots, g_{i+1})
+(-1)^{i+1}\sigma(g_1,\dots,g_{i}).
\end{split}
\end{align}
Set $Z^i(\Gamma, M)=\ker(d_i),$ $B^i(\Gamma, M)=d^{i-1} (C^{i-1}(\Gamma, M))$ and
$H^i(\Gamma, M):=Z^i(\Gamma, M)/B^i(\Gamma, M).$
For instance, a $1$-cocycle
 $\phi$ is a map from $\Gamma$ to $M$ such that
 \begin{equation}\label{1coc}
   \phi(g_2 g_1)=\phi(g_2).g_1+\phi(g_1) \qquad \qquad \text{for all $g_1, g_2 \in \Gamma$.}
 \end{equation}
Both $Z^i$ and $H^i$ are endowed with a Hecke action which we will not define but mention
because it plays an important role in an application below. Detailed expositions can be found in
\cite{ShimuraSur, D1}.

In Eichler cohomology, we apply this construction
with $M$ the space $P_{k-2}$ of polynomials of degree $\le k-2$.
The action $|_{2-k}$ of $\Gamma$ on $P_{k-2}$ or, more generally, on functions on $\mathfrak{H}$ is:
\begin{equation*}
(P|_{2-k} \gamma)(z):=P(\gamma z) j(\gamma, z)^{k-2}, \qquad z \in \HH, \gamma \in \Gamma,
\end{equation*}
where
$j\left ( \left ( \begin{smallmatrix} * & * \\ c & d \end{smallmatrix} \right ),
\tau \right ):=c\tau+d.$
With this notation, \eqref{cobound} is re-written as
\begin{equation}\label{cobound'}
F|_{2-k}S-F=r_f.
\end{equation}
We now consider the map $\sigma_f\colon \Gamma \to P_{k-2}$ defined by first setting
$\sigma_f(S)=r_f$ and $\sigma_f(T)=0$ and extending to $\Gamma$
according to \eqref{1coc}. This, in view of \eqref{ST}, gives a well-defined map because, with \eqref{cobound'},
\begin{equation}\label{welldef}
\sigma_f(S^2)=\sigma_f(S)|_{2-k}S+\sigma_f(S)=r_f|_{2-k}S+r_f=F|_{2-k}(S^2-S)+F|_{2-k}(S-1)=0,
\end{equation}
and likewise $\sigma_f((ST)^3)=0.$

Note that, despite the appearance of \eqref{cobound'}, $r_f$ is not a $1$-coboundary in $B^1(\Gamma, P_{k-2})$
because $F \not \in P_{k-2}$. It is only a coboundary in a larger space but, as we just saw, this fact suffices to
show that $r_f$ is $1$-cocycle in $Z^1(\Gamma, P_{k-2})$. This technique is used often in constructions of cocycles
and will reappear in the  sequel.

It is possible to express the value of $\sigma_f$ at every $\gamma \in \Gamma$ by a simple formula:
\begin{equation}
\label{sigmagamma}
\sigma_f(\gamma)(z)=\int_{\gamma^{-1} \infty}^{\infty}f(\tau) (\tau-z)^{k-2}dz.
\end{equation}

A fundamental fact is that Eichler cohomology parametrizes modular forms by means of the
{\it Eichler-Shimura isomorphism.} For general $f$ in the space $M_k$ of all weight $k$ modular forms for SL$_2(\mathbb{Z})$, 
we define $\sigma_f$ by 
\begin{equation*}
\sigma_f(\gamma)(z):=\int_{\gamma^{-1} \tau_0}^{\tau_0}f(\tau) (\tau-z)^{k-2}dz,
\end{equation*}
 where $\tau_0 \in \mathfrak{H}$ is fixed. If $\bar g$ is the function
obtained by conjugating the values of $g$, we define $r_{\bar f }$ by a similar formula involving integration of
antiholomorphic differentials. Then we have the following.
\begin{thm}{\rm (Eichler-Shimura isomorphism)} \label{ES} Let $\sigma$ be the map assigning to 
$(f, \bar g) \in M_k \oplus \overline{S_k} $ the $1$-cocycle $\sigma_f+\sigma_{\bar g}$ and let
$\pi$ be the natural projection of $Z^1(\Gamma, P_{k-2})$ onto $H^1(\Gamma, P_{k-2})$. Then $\pi \circ \sigma$ is a Hecke-equivariant isomorphism.
\end{thm}

From this viewpoint, the period polynomial of a cusp form $f$
can be redefined as the value at the involution $S$ of the
image of $f$ under the Eichler-Shimura map $\sigma$.

{\bf b. Critical values of $L$-functions}
A second interpretation of the period polynomial is as a
generating function of critical values of $L$-functions. 
As usual, we define the $L$-function of a modular form
$f(z)=\sum_{n=0}^{\infty}a_n e^{2 \pi i nz}$ by
$$L_f(s):=\sum_{n=1}^{\infty}\frac{a_n}{n^s} \quad\quad\text{ (for }\operatorname{Re}(s) \gg 0),$$ 
and the completed $L$-function by
$$\Lambda_f(s):=(2 \pi)^{-s} \Gamma(s) L_f(s).$$
The function
$\Lambda_f$ has a meromorphic continuation to the entire complex plane with possible
(simple) poles at $0$ and $k$, and it
 satisfies the functional equation (see, e.g., \cite{I}, Chapt. 7): 
\begin{equation}
\label{FE}
\Lambda_f(s)= i^k \Lambda_f(k-s).
\end{equation}
It further has an integral expression:
\begin{equation}\label{Mellin} \Lambda_f(s)=
\int_1^{\infty}(f(iv)-a(0))v^{s-1} dv+
i^k \int_1^{\infty}(f(iv)-a(0))v^{k-s-1} dv-
\frac{a(0)}{s}-\frac{a(0)i^k}{k-s}.
\end{equation}
When $f$ is a cusp form, $\Lambda_f(s)$ is entire and \eqref{Mellin} becomes
 the classical Mellin transform: 
\begin{equation}\label{Mellincusp}
\Lambda_f(s)=
\int_0^{\infty}f(iv)v^{s-1} dv.
\end{equation}

The values of the $L$-functions of Hecke eigencuspforms are of fundamental importance,
among other reasons, because they are, at least conjecturally, closely connected with
arithmetic and classical arithmetic questions. For example, a part of the
Birch--Swinnerton-Dyer conjecture implies that, if $L_f(1) \neq 0$, for a weight $2$ cusp
form of a certain type, then a specific polynomial Diophantine equation has at most finitely many
solutions. This is part of the order $0$ Birch--Swinnerton-Dyer conjecture and has
been proven in \cite{CoatesWiles}, \cite{GZ}, and \cite{Ko88, Ko89}.


Among the values of $L$-functions, the values at the integers within the critical strip
$0<\text{Re}(s)<k$ are called {\it critical} and were the first ones to be studied. 
Using the binomial theorem and \eqref{Mellincusp}, one can show that the period polynomial naturally encodes the critical $L$-values: 
\begin{equation}\label{CritPeriodP}
r_f(z)=-i \sum_{j=0}^{k-2} \binom{k-2}{j} (iz)^j \Lambda_f(j+1).
\end{equation}

\begin{rem}\label{defover}
{\rm The cohomological properties of $r_f$ discussed in Part a. can be translated to
analogous statements here. Notably, if we rewrite the equations
proving that $\sigma_f$ is well-defined (e.g. \eqref{welldef}) in terms
of  \eqref{CritPeriodP}, we are led to Manin's important
``Eichler-Shimura relations'' (Prop. 2.1 of \cite{M}).
A crucial implication of these relations, discussed in the next application, is that $Z^1(\Gamma, P_{k-2})$
(more precisely, the part of this space corresponding to the odd part
of the polynomials in $P_{k-2}$) is defined by a linear system of equations
with rational coefficients. 
}
\end{rem}

{\it Application of a. and b.} We will combine
the above two interpretations of $r_f$ above to illustrate the power of
the period polynomial with the following result.
\begin{thm}\label{ManinPeriodT} {\bf Manin's Periods Theorem} \cite{M}
Let $f$ be a normalized Hecke eigencuspform in $S_k$ and
let $K_f$ be the field obtained by adjoining to
$\Q$ the Fourier coefficients of $f$.
There exist $\omega_+(f)$, $\omega_-(f) \in \R$ such that
$$
\Lambda_f(s)/\omega_+(f), \quad \Lambda_f(w)/\omega_-(f) \in K_f
$$
for all $s,w$ with $1 \le s, w \le k-1$ and $s$ even, $w$ odd.
\end{thm}
Besides Manin's proof in \cite{M}, other proofs stressing different aspects
include those of Shokurov (\cite{Shok}, geometric methods on Kuga-Sato varieties),
Zagier (\cite{ZBonn}, using Rankin-Selberg method and Rankin-Cohen brackets),
Shimura (\cite{Shim}, by another variant of the Rankin-Selberg method),
the author and O'Sullivan (\cite{DOS}, by a variation of a method of
\cite{KoZ} which uses holomorphic projection and Cohen's kernel) etc.

{\it Sketch of proof of \ref{ManinPeriodT}} The Hecke eigencuspform $f$ generates
a one-dimensional eigenspace of $S_k$. By Th. \ref{ES}, this is mapped
isomorphically to a one-dimensional eigenspace of $H^1(\Gamma, P_{k-2})$ and, in fact,
the restriction of that map
to just the even (resp. odd) powers of the polynomial induces
an isomorphism too. It can also be proved that this map sends
$f$ to a one-dimensional eigenspace $A_f$ of the even (resp. odd)
part of $Z^1(\Gamma, P_{k-2})$ not just of
$H^1(\Gamma, P_{k-2})$.  By Remark \ref{defover}, $Z^1(\Gamma, P_{k-2})$ is defined
over $\Q$ and thus $A_f$ is defined over $K_f$. Since $\dim(A_f)=1$, this implies that
 there is a $c \in \C$ such that the even (resp. odd) part of $\sigma_f(S)=r_f$ equals
 $c P^+$ for an even polynomial $P^+ \in K[z]$ (resp. $c P^-$ for an odd polynomial $P^- \in K[z]$).
 With \eqref{CritPeriodP}, this implies that quotients of critical $L$-values of
 the same parity belong to $K_f$. \qed

\subsection{Period polynomials of non-cuspidal modular forms}
The definition \eqref{ppoly} no longer applies in the case that $f$ is not cuspidal because 
the integral may fail to
converge at the end points. However, it is possible to modify the definitions so that they
include general modular forms as well. To our knowledge, the first one to give a general
definition and systematically study it was Zagier in \cite{Z}. (Grosswald \cite{Gross}, starting 
from a different departure point, also worked with a similar object and proved a explicit expression for it.) 
In the case of general $f \in M_k$ the ``first definition'' \eqref{ppoly} was replaced
by
\begin{align*}
\widetilde r_f(z)&=&
\int_i^{\infty} (f(\tau)-a_0)(\tau-z)^{k-2}d\tau+
\int_0^i (f(\tau)-a_0\tau^{-k})(\tau-z)^{k-2}d\tau+ \nonumber
\\
\qquad \qquad \qquad &\,&\frac{a_0}{k-1}\left ( (z-i)^{k-1}+\frac{( 1+iz )^{k-1}}{z} \right ) \nonumber \\
&=&
\left(\int_{\infty}^{i} (f( \tau)-a_0)(\tau -z)^{k-2}d  \tau 
+ \frac{a_0}{k-1}(i-z)^{k-1} \right)\Bigg|_{2-k}(S-1).
\end{align*}
The ``second definition'' \eqref{CritPeriodP}
was replaced by
\begin{equation}\label{CritPeriodPE}
\widetilde r_f(z)=
-i \sum_{j=0}^{k-2} \binom{k-2}{j} (iz)^j \Lambda_f(j+1)+
\frac{a_0}{k-1}\left ( z^{k-1}+z^{-1} \right ).
\end{equation}
This extended definition is then used in \cite{Z} to state and prove an
expression of a striking generating function involving period polynomials
over a basis of $M_k$ as a quotient of
products of values of the classical Jacobi theta function.

A difference from the case of cusp forms is that $\widetilde r_f$ is not in $P_{k-2}$ when $f$ is not cuspidal.
Recently,  it was shown in \cite{Br} that it is possible to define the period polynomial of
Eisenstein series so that it stays within $P_{k-2}$. Set
\begin{equation}\label{ppolyEBr}
r_f(z)=  \left(\int_{\infty}^z (f(w)-a_0)(w-z)^{k-2}dw
+a_0 \int_0^z (w-z)^{k-2}dw \right)\Bigg|_{2-k}(S-1).
\end{equation}
This definition was made in a vastly general context by Brown in \cite{Br} which
included general iterated Shimura integrals and originated from an
integral at a tangential base point at infinity. In this more general setting it is proved that
$r_f \in P_{k-2}$ and that it induces a $1$-cocycle. Brown's extension of the
period polynomial was also motivated by important applications. For example, he used it
to express non-critical values in terms of multiple
modular values.

A way to compare this definition with that of \cite{Z} is to consider the
explicit form of the Eichler cocycle it induces. Set
$$v_f(z):=\int_{\infty}^z (f(w)-a_0)(w-z)^{k-2}dw
+ \frac{a_0}{k-1}z^{k-1}.$$
Then, with the definition of $d^0$ in \eqref{differential} we define $\sigma_f:=d^0 v_f$. 
 It is clear that $r_f=\sigma_f(S)$. On the other hand, the associated cocycle of $\widetilde r_f$ is given by
$\widetilde \sigma_f:=d^0 \widetilde{v}_f$,
where
$$\widetilde{v}_f(z):=\int_{\infty}^i (f(w)-a_0)(w-z)^{k-2}dw
+ \frac{a_0}{k-1}(i-z)^{k-1}$$
and where $d^0$ is defined by the same formula as \eqref{differential}
but its domain is enlarged to $C^0(\Gamma, \C (z)).$

As mentioned above, in \cite{Br} it is proven in more
general form that
$\sigma_f$ takes values in $P_{k-2}$. This can also be seen by the identity
shown in \cite{DR} (eq. (8)):
\begin{align}\label{Eichlercocycle}
\begin{split}
\sigma_f(\gamma)
&=\int_{\gamma^{-1} i}^{i} f(w)(w- z)^{k-2}dw\\
&+
\left(\int_{\infty}^i (f(w)-a_0)(w-z)^{k-2}dw
+a_0 \int_0^i (w-z)^{k-2}dw \right)\Bigg|_{2-k}(\gamma-1).
\end{split}
\end{align}
From this it is also clear that $\sigma_f$ is
``canonical'' in the sense that it belongs to the same cohomology class
as the image of $f$ under the Eichler-Shimura isomorphism (Th. \ref{ES}).
With this definition of the cocycle $r_f$ we then have
\begin{equation}\label{LderBr}
r_f(z)=-i \sum_{j=0}^{k-2} \binom{k-2}{j} (iz)^j \Lambda_f(j+1)
\end{equation}
(see the proof of Lemma 7.1 in \cite{Br})
Consider the case when $f=E_k$, the Eisenstein series
 $$E_k(\tau)=-\frac{B_k}{2k}+\sum_{n=1}^{\infty} \sigma_{k-1}(n)e^{2 \pi i n \tau},$$
where $B_a$ is the $a$-th Bernoulli number. Then 
Lemma 7.1 of \cite{Br} shows 
that the above equation takes the form
\begin{equation}
\label{Bernoulli}
r_{E_k}(z)=- \frac{(k-2)!}{2}\sum_{j=0}^{k/2-2} \frac{B_{2j+2}}{(2j+2)!} \frac{B_{k-2j-2}}{(k-2j-2)!}
z^{2j+1}+\frac{(k-2)!}{2} \frac{\zeta(k-1)}{(2 \pi i)^{k-1}}  (1-z^{k-2}),
\end{equation}
and, by Proposition of pg. 453 of \cite{Z} (or \eqref{CritPeriodPE}), we have 
\begin{equation}
\label{Bernoullitilde}
\widetilde r_{E_k}(z)= - \frac{(k-2)!}{2}\sum_{j=-1}^{k/2-1} \frac{B_{2j+2}}{(2j+2)!} \frac{B_{k-2j-2}}{(k-2j-2)!}
z^{2j+1}+\frac{(k-2)!}{2} \frac{\zeta(k-1)}{(2 \pi i)^{k-1}} (1-z^{k-2}).
\end{equation}

\subsection{Zeros of period polynomials}
\label{ZerosPeriodPolysSubSection}
Having argued the case for the conceptual importance of period polynomial and for its usefulness
due to its structure, it becomes clear that it is of interest to examine it for its own sake, as a polynomial.
We review work on its zeros as a polynomial.

For $k \in 2 \mathbb{N}$, M. R. Murty, C. Smyth and R. Wang \cite{MSW} studied the {\it Ramanujan polynomial}
\begin{equation}
\label{RamanPoly}
\sum_{j=0}^{k/2} \frac{B_{2j}}{(2j)!} \frac{B_{k-2j}}{(k-2j)!}
z^{2j}
\end{equation}
and proved the following result. 
\begin{theorem}[M. R. Murty, C. Smyth, R. Wang \cite{MSW}]\label{msw}
All non-real zeros of the Ramanujan polynomial are on the unit circle. 
\end{theorem}
With \eqref{Bernoullitilde},  the Ramanujan polynomial equals
the odd part of $-2z \widetilde r_{E_k}(z)/(k-2)!$. 

Because of \eqref{welldef}, the circle is a natural ``line of symmetry'' for the period polynomials,
and therefore results such this can be thought of as a ``Riemann Hypothesis'' for period polynomials.
This viewpoint was adopted in \cite{ORS} where similar statements are connected to
Manin's theory of ``zeta polynomials'' $Z_f(s)$. These are versions of the period polynomials
that send the unit circle to $\operatorname{Re}(s)=\frac12$ and satisfy
the functional equation $Z_f(1-s)=\pm Z_f(s)$.

In \cite{LS}, Lal\'\i n and Smyth studied the zeroes of the ``Ramanujan polynomials''
\begin{equation*}
\mathcal R_k(z):=\sum_{j=0}^{k/2} \frac{B_{2j}}{(2j)!} \frac{B_{k-2j}}{(k-2j)!}
z^{2j}+\frac{\zeta(k-1)}{(2\pi i)^{k-1}}(z^{k-1}-z).
\end{equation*}
\begin{theorem}[Lal\'\i n and Smyth \cite{LS}]
\label{EisPolyTheorem}
For each $k \in 2 \mathbb{N}$, the zeroes of $\mathcal R_k$ all lie on the unit circle.
\end{theorem}
From a modular perspective, our interest in these polynomials is that, by \eqref{Bernoullitilde} 
they are the full $-2z \widetilde r_{E_k}(z)/(k-2)!$. In the sequel, we shall be concerned with analogues 
of Th.~\ref{msw} for
derivatives of Eisenstein $L$-series. In particular, it will serve as the main motivation for 
our first steps towards understanding our broader conjectures for entire modular forms spaces.

An interesting recent interpretation of the period polynomial of Eisenstein
series and of their zeros in view of Ramanujan's ``formula'' for $\zeta(2m+1)$ is discussed in
\cite{BS}. In the same paper a question is raised about a variation of the Ramanujan
polynomial:
\begin{equation}\label{p_m}
p_m(z)=\frac{\zeta(2m+1)}{2}(1-z^{2m})-\frac{(2 \pi i)^{2m+1}}{2}\sum_{n=1}^{m}
\frac{B_{2n}}{(2n)!}\frac{B_{2m-2n+2}}{(2m-2n+2)!} z^{2n-1}.
\end{equation}
Note that $p_{k/2-1}(z)=(2 \pi i)^{k-1} r_{E_k}(z)/(k-2)!.$ In Remark 7.4 of \cite{BS}, 
it is asked whether $p_m$ and $p_m^-(z)/z$ (where $p_m^-$ is the odd part of $p_m$) are unimodular. In \cite{DR}, we proved
 the second part of this conjecture.

Analogous results have been proved
for cusp forms. For example, Conrey, Farmer, and Imamo{\=g}lu \cite{CFI}
have proved that, apart from five ``trivial'' real zeroes, all zeroes of the odd part of the
period polynomial of a cusp form lie on the unit circle.
\begin{theorem}[Conrey, Farmer, and Imamo{\=g}lu]
\label{OddPolyTheorem}
If $f$ is a cuspidal Hecke eigenform 
on $\mathrm{SL}_2(\ZZ)$, then the odd part of $r_f$ has zeroes at $0,\pm\frac12,\pm2$. The remainder of the zeros lie on the unit circle.
\end{theorem}

A similar picture exists for the full period polynomials $r_f$ for Hecke eigencuspforms $f$. However, in this case, there are no trivial zeros, and all zeros of $r_f$ lie on the unit circle. This
is summarized in the following result, shown by El-Guindy and Raji in \cite{ElG-R} for level $1$ and for general level $N$ by Jin, Ma, Ono, and Soundararajan in \cite{JMOS}.
\begin{theorem}[El-Guindy and Raji and Jin, Ma, Ono, and Soundararajan]
\label{FullPolyTheorem}
If $f$ is a Hecke eigencuspforms on $\Gamma_0(N)$ for any $N$, then all zeroes of the period polynomial $r_f$ lie on the unit circle.
\end{theorem}
\begin{remark}
{\rm Explicit approximations for the exact locations of the zeroes were given in \cite{JMOS}. }
\end{remark}

The proofs of the above results are based on the origin of the
period polynomial as a cocycle. In particular, the behaviour of $r_f$ under the action of the involution $S$
imposes a special structure on the polynomial (it is a {\it self-inversive polynomial}). This allows for a more
convenient investigation of the location of the zeros thanks to the following result.

\begin{lemma}[Th. 2.2 of \cite{ElG-R}]
If $h(z)$ is a polynomial of degree $n$ with all zeros inside the unit disk $|z|\leq1$, then for any $d\geq n$ and $\lambda$ on the unit circle, 
the polynomial
\begin{equation}
\label{ElGRTheorem}
z^{d-n}h(z)+\lambda z^n\overline{h}(1/z)
\end{equation}
has all its zeroes on the unit circle, provided that it is not identically zero.
\end{lemma}
Statements of this type have a long history which can be traced back to Hermite (see the Addendum of \cite{LS} for an account) but, in this form, the proposition has been proved in \cite{ElG-R}.

In this way, the problem of locating the zeros of $r_f$ is reduced to locating the zeros of the polynomial $h$ associated to $r_f$ through 
Lemma \ref{ElGRTheorem}. This is achieved by bounds and monotonicity statements for values of $L$-functions 
appearing in the coefficients of $h$. 
The results cited above are proved by using different such bounds and monotonicity statements.

\section{``Period polynomial'' for derivatives of $L$-functions}\label{Der}
Beilinson's conjecture, part of which is stated in Conj. \ref{Beil}, pertains to 
values of derivatives of L-functions and, as mentioned in the introduction, very 
little is known about the case of order greater than $0$. This has motivated many
approaches to the study of values of derivatives. We will outline one, due to Goldfeld and 
the first author (see \cite{G, D1, D} and the later works by them and their collaborators: \cite{D2, BCD, DNS}) that
incorporates these values into the Eichler cohomology setting. In the cited papers,
only the cuspidal case was studied but here we will describe the general case as that
was described in \cite{DR}. 

\subsection{First derivatives}\label{Firstderivatives}
 We first recall the Dedekind eta function $$\eta(\tau):=e^{\frac{2 \pi i \tau}{24}} \prod_{n=0}^{\infty} (1- e^{2 \pi i n \tau})$$
and then set
$u(\tau):=2\log(\eta(\tau))$.
For each $\gamma \in \Gamma$, this function satisfies
\begin{equation}\label{v}
u(\gamma \tau)=u(\tau)+\log (j(\gamma, \tau))+c_{\gamma}
\end{equation}
for some $c_{\gamma} \in \C$. In particular, $c_S=-\frac{\pi i}{2}$.

With the definition \eqref{differential} we set $\sigma_f:=d^1v_f$ where 
\begin{align*}
v_f(\gamma)
&
:=\int_{\infty}^z (f(w)-a_0)(w-z)^{k-2}
\left ( u(\gamma w)-u(w) \right )dw
\\
&+
a_0 \int_i^z (w-z)^{k-2}\left ( u(\gamma w)-u(w) \right )dw.
\end{align*}
It can be proved that, 
although $v_f$ is a cochain that takes values in the space $\mathcal{O}$ of
holomorphic functions on the upper-half plane, $\sigma_f$ takes values in
the much smaller space of polynomials of degree $\le k-2$. Since, further, 
it is in the image of the differential map $d^1$ , we deduce:
\begin{lem}\label{Lemfirst} (Lemma 3.3 of \cite{DR})The map $\sigma_f$
is a a $2$-cocycle in $P_{k-2}$.
\end{lem}
As mentioned above, this construction extends the corresponding one or
cusp forms given in \cite{D}. This is the content of the following
proposition which, further, expresses $\sigma_f$ in a way which make the
analogy with the standard polynomial \eqref{sigmagamma} more transparent. 
\begin{prop} Let $f$ be a cusp form of weight $k$ for $\Gamma$.  Then
\begin{align*}
\begin{split}
\sigma_f(\gamma_1, \gamma_2)
&=\int_{\gamma^{-1}_1 \infty}^{\infty} f(w)(w-z)^{k-2}(u(\gamma_2 w)-u(w))dw\\
&=\int_{\infty}^{\gamma_1 \infty} f(w)(w-z)^{k-2}(u(\gamma_2 w)-u(w))dw\Big |_{2-k}\gamma_1.
\end{split}
\end{align*}
\end{prop}
The connection with values of derivatives of $L$-functions
is given by 
\begin{prop}\label{1stder} (Prop. 3.5 of \cite{DR}) Set $$P(z)=
\sum_{n=0}^{k-2} \binom{k-2}{n} \frac{i^{1-n}}{(n+1)^{2}}z^{k-2-n}.$$ Then
$$\sigma_f(S, S)=
-\sum_{n=0}^{k-2} \binom{k-2}{n} i^{1-n}\Lambda_f'(n+1)z^{k-2-n}
+a(0) (P|_{2-k}(1+S))(z).$$
\end{prop}
The proposition is stated in general, but, for cuspidal $f$, the analogy to \eqref{LderBr} is obvious. 

\subsection{Zeros of ``period polynomials'' for $L$-derivatives}\label{0s}
In light of the analogy with the standard period polynomial, 
it is natural to ask whether similar patterns in the distribution of 
zeros occur in ``period polynomials'' for derivatives of $L$-functions. 
Inspired by the behavior exhibited by ordinary period 
polynomials as described in Theorems~\ref{OddPolyTheorem} and \ref{FullPolyTheorem},
the authors searched for similar properties for polynomials  built from 
$L$-derivatives in \cite{DR}. There, the analogous period polynomials were defined to be
polynomials 
\[
Q_f(z):=\sum_{n=0}^{k-2} \binom{k-2}{n} i^{1-n}\Lambda'_f(n+1)z^{k-2-n},
\]
in direct analogy with \eqref{LfPeriodP} and the following conjecture was formulated.
\begin{conjecture}\label{Conj0first} 
For any Hecke eigenform of weight $k$ on $\operatorname{SL}_2(\mathbb Z)$ the polynomial $Q_f(z)$
has all its zeros on the unit circle. 
Moreover, the odd part of $Q_f(z)$ has all of its zeros on the unit circle, except for trivial
zeros at $0$ and $\pm a,\pm1/a$ for some real number $a$. 
\end{conjecture}
The evidence for this conjecture was both theoretical and experimental. The former  
was provided by 
our proof of the second part of Conj. \ref{Conj0first} in the case of Eisenstein series. This is, at the same time, the analogue
of the main result of \cite{MSW} on period polynomials of Eisenstein series. 
\begin{thm}\label{1stderE*} (\cite{DR}) If $4|k$, all non-zero zeroes of the odd part of $Q_{E_k}$
lie on the unit circle. 
\end{thm}
As in the case of the standard period polynomial, the pivot of the proof is the cohomological origin of the ``period polynomial'' $Q_{E_k}$ 
which allows us to study it as a self-inversive polynomial. On the other hand, we were then able to use more general theorems about 
locations of zeros (Enestr\"om-Kakeya Theorem \cite{En,Ka}). This is because our construction parallels Brown's version of the period polynomial of Eisenstein series ($r_{E_k}$) and not that of \cite{Z} ($\tilde{r}_{E_k}$) which, with its two extra terms, takes us away from
the coefficient module of polynomials.

At first glance, since in Th. \ref{msw} there are further real roots (in addition to $0$), the conclusion of Th. \ref{1stderE*} appears to not be analogous with its counterpart Th. \ref{msw}. The reason
for this is that, whereas the subject of Th. \ref{msw} is Zagier's version of the period polynomial of Eisenstein series, the subject of Th. \ref{1stderE*} is a polynomial which extends Brown's version of the period polynomial of Eisenstein series. The analogue of Th. 
\ref{msw} for Brown's version of the period polynomial of Eisenstein series was stated as a question in \cite{BS} and has the same conclusion
as Th. \ref{1stderE*} (as shown in \cite{DR}).
 
\noindent
{\bf Question:} {\it Are all non-zero zeroes of the odd part of 
$$p_{k/2-1}(z)=\frac{(2 \pi i)^{k-1}}{(k-2)!} r_{E_k}(z)$$ on the unit circle?  }

\medskip 

The experimental evidence for the truth of Conj. \ref{Conj0first} is also very convincing and will be outlined along the respective discussion of higher derivatives.

We end this section by noting that it would be very interesting to interpret the role of the number $a$ in the statement of the conjecture, 
and in particular to find an explanation for them as ``trivial zeros'', as was the case for the zeros with $a=2$ in Th.~\ref{OddPolyTheorem}.

\subsection{The case of higher derivatives}
An advantage of the approach on derivatives of $L$-function discussed here is that
it includes in a natural way higher derivatives about which, as mentioned earlier, very little
is known. Therefore, any progress by this method in the case of first derivatives might lead
to insights for higher derivatives as well. 

The cohomological tool enabling to extend the constructions of Section \ref{Firstderivatives} to 
higher derivatives is {\it cup products}. This, in the case we need it, is defined as a map
$$\cup\colon C^1(\Gamma, \mathcal O) \otimes
C^m(\Gamma, \mathcal O) \to C^{m+1}(\Gamma, \mathcal O) $$
given by
$$\left ( \phi_1 \cup \phi_2 \right )(\gamma_1, \gamma_2, \dots, \gamma_{m+1}):=\phi_1(\gamma_1) \left ( \phi(\gamma_2, \dots, \gamma_{m+1})|_0 \gamma_1\right ).$$
A crucial property that that cup products of cocycles are cocycles.
For $\phi_i \in C^1(\Gamma, \mathcal{O})$, we set:
$$\phi_1 \cup \dots \cup \phi_n:=\phi_1 \cup \left(\phi_2 \cup \left( \dots ( \phi_{n-1} \cup \phi_n ) \dots \right)\right) \in C^n(\Gamma, \mathcal O).$$

If $v$ is the $1$-cocycle given by
$\gamma \to u|_0(\gamma-1)$ (with $u$ as in Section \ref{Firstderivatives}), we 
 set, for $n \in \mathbb N$,
$$V_n:=v \cup v \cup \dots \cup v \qquad \text{($n$ times)}.$$
As mentioned above, this will be a $n$-cocycle.

Let $v_f \in C^n(\Gamma, \mathcal{O})$ be given by
\begin{align*}
v_f(\gamma_1, \dots, \gamma_n)
&=\int_{\infty}^z (f(w)-a_0)(w-z)^{k-2} V_n(\gamma_1, \dots, \gamma_n)(w)dw\\
&+ a_0\int_{i}^z (w-z)^{k-2} V_n(\gamma_1, \dots, \gamma_n)(w)dw.
\end{align*}
Setting $\sigma_f:=d^n v_f$, we arrive at the following analogue of Lemma~\ref{Lemfirst} for higher cocycles.
\begin{lem}\label{mthderiv} (Lemma 3.6 of \cite{DR}) The map $\sigma_f$
takes values in $P_{k-2}$ and thus gives an $(n+1)$-cocycle in $P_{k-2}$.
\end{lem}
Finally, the analogue of Prop. \ref{1stder} is
\begin{prop}\label{higherderivative} (Prop. 3.7 of \cite{DR}) For each $m \in \mathbb{N}$, set $$P(z)=
\sum_{n=0}^{k-2} \binom{k-2}{n} \frac{i^{1-n}}{(-n-1)^{m+1}}z^{k-2-n}.$$ Then
\begin{equation}\label{Q_f}
(-1)^{m}\sigma_f(S, \dots S)=
\sum_{n=0}^{k-2} \binom{k-2}{n} i^{1-n}\Lambda_f^{(m)}(n+1)z^{k-2-n}
-a(0) m! (P|_{2-k}(1+(-1)^{m+1}S))(z),
\end{equation}
where $\sigma_f$ has $m+1$ arguments.
\end{prop}
This proposition led us to formulate Conj. \ref{Conj0} as the general version of Conj. \ref{Conj0first}.
\begin{conjecture}
\label{FinalConjecture}
For any Hecke eigenform of weight $k$ on $\operatorname{SL}_2(\mathbb Z)$, and for each $m\in\mathbb Z_{\geq0}$, the polynomial
\[
Q_f(z):=\sum_{n=0}^{k-2} \binom{k-2}{n} i^{1-n}\Lambda_f^{(m)}(n+1)z^{k-2-n}
\]
has all its zeros on the unit circle.
Moreover, the odd part
\[\sum_{\substack{n=1 \\ n \, \, \text{odd}}}
^{k-3}\binom{k-2}{n} i^{1-n}\Lambda_f^{(m)}(n+1)z^{k-2-n}
\]
has all of its zeros on the unit circle, except for zeros at $0$ and $\pm a,\pm1/a$ for some $a\in\mathbb R$. 
\end{conjecture}

We were able to prove the Eisenstein series case of the second part of this conjecture too, 
but we had to truncate the ``lower order'' terms from the $\Lambda^{(m)}_f$ appearing in $Q_f$. The precise
construction is slightly complicated but the essence of the theorem is entirely analogous to 
Th. \ref{1stderE*} (see Th. 4.2 of \cite{DR}).

The experimental evidence for Conj.~\ref{FinalConjecture} in the case of both the first and the higher derivatives was based on computer
search. In particular, the authors used SAGE to check that the norms of all zeroes of all full period polynomials with $m\leq3$ and $k\leq 50$ were within $10^{-10}$ of $1$. The structure of the second part of the conjecture was made on the basis of similar computational experiments.


\begin{thebibliography}{99}

 \bibitem{BCD}
R. Bruggeman, Y. Choie, and N. Diamantis, \emph{Holomorphic automorphic forms and cohomology}, Memoirs of the AMS (to appear)
arXiv:1404.6718.

 \bibitem{BS}
B. Berndt, A. Straub \emph{Ramanujan's Formula for $\zeta(2n+1)$}, arXiv:1701.02964v1.

 \bibitem{Br}
F. Brown, \emph{Multiple Modular Values and the relative completion of the fundamental group of $M_{1,1}$}, preprint,
 arXiv:1407.5167.

\bibitem{CoatesWiles} J. Coates, A. Wiles, \emph{On the conjecture of Birch and Swinnerton-Dyer},
Inventiones Math. {\bf 39} (1977) 223--251.

\bibitem{Cohn} A. Cohn, \emph{\"Uber die Anzahl der Wurzeln einer algebraischen Gleichung in einem Kreise}, Math. Zeit. {\bf 14}
(1922), 110--148.

 \bibitem{CFI} J.B. Conrey, D.W. Farmer, and {\"O}. Imamo{\=g}lu, {\it The nontrivial zeros of period polynomials of modular forms lie on
the unit circle}, Int. Math. Res. Not. no. 20, 4758--4771  (2013).

\bibitem{De}
P. Deligne \emph{Valeurs de Fonctions L et p\'eriodes d'int\'egrales},
Proceedings of Symposia in Pure Mathematics {\bf 33}, 313-346 (1979).

\bibitem{D1}
N. Diamantis,\emph{
Special values of higher derivatives of $L$-functions},
Forum Math. \textbf{11} no. 1, 229--252 (1999).

\bibitem{D}
N. Diamantis,\emph{
Hecke Operators and Derivatives of $L$-Functions},
Compositio Math. \textbf{125} no. 1, 39--54 (2001).


\bibitem{D2}
N. Diamantis,\emph{The geometry of certain cocycles associated to derivatives of $L$-functions},
Forum Math. \textbf{17} no. 5, 739--752 (2005).


\bibitem{DNS}
N. Diamantis, M. Neururer, F. Str\"omberg, \emph{A correspondence of modular forms and applications to values of $L$-series},
Res. number theory (2015) 1: 27.

\bibitem{DOS}
N. Diamantis, C. O'Sullivan, \emph{Kernels of $L$-functions of cusp forms} Math. Ann. 346 (2010), no. 4, 897--929.

\bibitem{DR}
N. Diamantis, L. Rolen, \emph{Eichler cohomology and zeros of polynomials associated to derivatives of $L$-functions},
arXiv:1704.02667.



\bibitem{DIT}
W. Duke, {\"O}. Imamo{\=g}lu, \'A. T\'oth,  \emph{Rational period functions and cycle integrals},
Abhandlungen aus dem Mathematischen Seminar der Universit\"at Hamburg
 {\bf 80} (2) (2010) 255--264.

\bibitem{EichlerVer} M. Eichler, \emph{Eine Verallgemeinerung der Abelsche Integrale} Math. Z. {\bf 67} (1957) 267--€"298.

\bibitem{ElG-R}
A. El-Guindy, W. Raji, \emph{Unimodularity of zeros of
period polynomials of Hecke eigenforms}, Bull. Lond. Math. Soc. \textbf{46} no. 3, 528--536 (2014).


\bibitem{En} G. Enestr\"om, \emph{Ramarque sur un th\'eor\`eme relatif aux recines de l'equation $a_nx_n+\cdots+a_0 = 0$ o\`u tous les coefficients sont r\'eels et positifs}, T\^{o}hoku Math. J. \textbf{18}, 34--36 (1920), translation of a Swedish article in Ofversigt of Konogl. Vertenskaps Akademiens F\"orhandlingar \textbf{50}, 405--415 (1893).




\bibitem{G}
D. Goldfeld, \emph{Special values of derivatives of $L$-functions},  Number theory
(Halifax, NS, 1994), 159--173, CMS Conf. Proc., \textbf{15}, Amer. Math. Soc., Providence, RI (1995).

\bibitem{GZ}
B. Gross and D. Zagier \emph{Heegner points and derivative of $L$-series}
Invent. Math. {\bf 85}, 225--320 (1986).

\bibitem{Gross}
E. Grosswald, \emph{Die Werte der Riemannschen Zeta-funktion an ungeraden Argumentstellen}
Nachr. Akad. Wiss. G\"ottingen (1970), 9--13.


\bibitem{I} H. Iwaniec, {\it Topics in classical automorphic forms}
Graduate Studies in Mathematics, Vol. 17, AMS, 1991.

\bibitem{JMOS} S. Jin, W. Ma, K. Ono, and K. Soundararajan, {\it The Riemann Hypothesis
for period polynomials of modular forms}, Proc. Natl. Acad. of Sci. U.S.A. \textbf{113} no. 10, 2603--2608 (2016).

\bibitem{Ka} S. Kakeya, \emph{On the limits of the roots of an algebraic equation
with positive coefficients}, Tohoku Math. J. \textbf{2}, 140--142 (1912-1913).

\bibitem{KoZ} W. Kohnen and D. Zagier, {\it Modular forms with rational periods
in Modular Forms}, R.A. Rankin (ed.), Ellis Horwood, Chichester 197--249 (1984).

\bibitem{KZ} M. Kontsevich and D. Zagier, {\it Periods}, Mathematics unlimited --2001 and beyond, 771--808, Springer, Berlin, 2001.

\bibitem{Ko88} V. Kolyvagin. {\it Finiteness of $E(\mathbb{Q})$ and $\Sha(E, \mathbb{Q})$ for a subclass of Weil curves.}
Izv. Akad. Nauk SSSR Ser. Mat. {\bf 52} (1988), no. 3, 522--540.

\bibitem{Ko89} V. Kolyvagin. {\it The Mordell-Weil and Shafarevich-Tate groups for Weil elliptic
curves}, Izv. Akad. Nauk SSSR Ser. Mat. {\bf 52} (1988), no. 6, 1154--€"1180.


\bibitem{LS} M. Lal\'\i n and C. Smyth,
{\it Unimodularity of zeros of self-inversive polynomials} Acta Math. Hungar. {\bf 138}
(2013), no. 1-2, 85--101. Addendum, Acta Math. Hungar. {\bf 147} (2015), no. 1, 255--257.

\bibitem{M} Y. T. Manin, \emph{Periods of parabolic points and $p$-adic Hecke series,} Math. Sb.,
371--393 (1973).


\bibitem{MSW} M. Murty, C. Smyth, and R. Wang, {\it Zeros of Ramanujan polynomials}
J. of the Ramanujan Math. Soc. \textbf{26}, 107--125 (2011).

\bibitem{ORS} K. Ono, L. Rolen, and F. Sprung, {\it Zeta-polynomials for modular form periods}, Adv. Math. \textbf{306}, 328--343 (2017).


\bibitem{NIST} F. Olver, D. Lozier, R. Boisvert, and
C. Clark, \textsl{NIST handbook of mathematical
functions} U.S.
Department of Commerce, National Institute of Standards
and Technology,
Washington, DC; Cambridge University Press, Cambridge, 2010.

\bibitem{P} H. Poincar\'e,  {\it Sur les invariants arithm\`etiques} J. f\"ur Reine u. Ang. Math. {\bf 129}, 89--150 (1905).

\bibitem{Shok} V. Shokurov,  {\it Shimura integrals of cusp forms}, Izv. Akad. Nauk SSSR Ser. Mat. {\bf 44} (1980), no. 3, 670--718, 720.

\bibitem{ShimuraSur} G. Shimura \emph{Sur les int\'egrales attach\'ees aux formes automorphes} J. math. Soc. Japan {\bf 11}
(1959) 291--311

\bibitem{Shim} G. Shimura,
{\it The special values of the zeta functions associated with cusp forms}
Comm. Pure Appl. Math. {\bf 29} (1976), no. 6, 783--804.


\bibitem{Z} D. Zagier, \emph{ Periods of modular forms and Jacobi theta functions},
Invent. Math. \textbf{104} no. 3 449--465 (1991).

\bibitem{ZBonn} D. Zagier, \emph{Modular forms whose Fourier coefficients involve zeta-functions of quadratic
fields} Modular functions of one variable, VI (Proc. Second Internat. Conf., Univ.
Bonn, Bonn, 1976), LNM 627, Springer, Berlin, 1977, 105--169.

\end{thebibliography}
\end{document}